\documentclass[12pt, titlepage]{amsart}

\usepackage{ae} 
\usepackage[T1]{fontenc}
\usepackage[cp1250]{inputenc}
\usepackage{amsmath}
\usepackage{amssymb, amsfonts,amscd,verbatim}
\usepackage[dvips]{graphicx}
\usepackage{latexsym}
\usepackage{indentfirst}
\usepackage{latexsym}

\usepackage{graphicx}
\usepackage{verbatim}   
\usepackage{color}      
\usepackage{subfigure}  

\usepackage{amsmath}    

\theoremstyle{plain}
\newtheorem{Prop}{Proposition}[section]
\newtheorem{Thm}[Prop]{Theorem}
\newtheorem{Cor}[Prop]{Corollary}

\newtheorem{Lem}[Prop]{Lemma}

\theoremstyle{definition}
\newtheorem{Def}[Prop]{Definition}

\theoremstyle{remark}

\newtheorem{Problem}[Prop]{\bf Problem}

\newcommand{\edim}{{\mathrm{edim}}}

\def\int{\mathop{\roman{int}}}

\def\1{^{-1}}

\def\Z{{\mathbf Z}}

\def\dim{\text{dim}}

\def\dim{\mathrm{dim}}

\def\Map{\mathrm{Map}}

\newcommand{\rp}{{\Bbb R}P}
\def\dokaz{{\bf Proof. }}
\def\enddokaz{\hfill $\blacksquare$}

\numberwithin{equation}{section}

\input pstricks.tex
\input xy
\xyoption{all}

\begin{document}
\title[
Maps to the projective plane
]%
   {Maps to the projective plane}
                  
\author{Jerzy Dydak}                  
\address{ Department of Mathematics\\
 University of Tennessee\\
 Knoxville, TN 37996-1300 }                  
\email{ dydak@math.utk.edu }                     

\author{Michael  Levin} 
\address{ Department of Mathematics\\
Ben Gurion University of the Negev\\
P.O.B. 653\\
Be'er Sheva 84105, Israel}             
\email{mlevine@math.bgu.ac.il } 
\date{ \today
}


\thanks{This research was supported
by Grant No.  2004047  from the United States-Israel Binational Science
Foundation (BSF),  Jerusalem, Israel.}
\thanks{The first-named author was partially supported
by the Center for Advanced Studies in Mathematics
at Ben Gurion University of the Negev (Beer-Sheva, Israel).}

\begin{abstract}

We prove the projective plane $\rp^2$ is an absolute extensor of a finite-dimensional
metric space $X$ if and only if the cohomological dimension mod $2$
of $X$ does not exceed $1$. This solves one of the remaining difficult
problems (posed by A.N.Dranishnikov) in extension theory.
One of the main tools is the computation of the fundamental
group of the function space $\Map(\rp^n,\rp^{n+1})$ (based at inclusion)
as being isomorphic to either $\Z_4$ or $\Z_2\oplus\Z_2$ for $n\ge 1$.
Double surgery and the above fact yield the proof.
\end{abstract}
\maketitle

\medskip
\medskip
\tableofcontents

\section{Introduction}

The basic relation studied by Extension Theory
is that of a CW complex $K$ being an {\it absolute extensor} of a metric space $X$.
It means that every map $f\colon A \to K$, $A$ closed in $X$, extends over $X$.
There are three existing notations for that:
\begin{enumerate}
\item $K\in AE(X)$.
\item $X\tau K$.
\item $\edim X \leq K$.
\end{enumerate}
In addition, we will use notation $X\tau (L\to K)$, where $L$ is a subcomplex
of $K$. It means that every map $f\colon A \to L$, $A$ closed in $X$, extends over $X$
with values in $K$. See \cite{CDSVV} for more information about that relation.
\par In case of basic CW complexes the relation $X\tau K$
is equivalent to classical concepts in Dimension Theory (see \cite{dr3} for more details):
\begin{itemize}
\item[a.] $X\tau S^n$, $S^n$ being the $n$-sphere, is equivalent to 
covering dimension $\dim(X)$ being at most $n$.
\item[b.] $X\tau K(G,n)$, $K(G,n)$ being the Eilenberg-MacLane complex,
is equivalent to the cohomological dimension $\dim_G(X)$ being at most $n$.
\end{itemize}
The leading theme in Extension Theory is the effort to relate
$X\tau K$ to a set of conditions $X\tau K(G_n,n)$, where $G_n$ depends
on $K$. In that vain
Dranishnikov \cite{dr2} proved the following
  important  theorems connecting  extensional and cohomological
  dimensions.

  \begin{Thm}
  \label{t-1}
  Let $K$  be a CW-complex and let a compactum $X$ be such that
  $\edim X \leq K$. Then
  $\dim_{H_n(K)}(X) \leq n$ for every $n>0$.
  \end{Thm}

   \begin{Thm}
  \label{t-2}
  Let $K$  be a simply connected CW-complex and let a compactum
  $X$ be finite
  dimensional.
  If $\dim_{H_n(K)}(X) \leq n$ for every $n>0$, then
  $\edim X \leq K$.
 \end{Thm}
 Both Theorems \ref{t-1} and \ref{t-2} were subsequently generalized
 for metric spaces $X$ in \cite{Dy1}.
 The requirement in Theorem \ref{t-2} that $X$ is finite dimensional
  cannot be omitted. To show that  take
 the famous infinite-dimensional compactum $X$ of
  Dranishnikov
   with $\dim_\Z (X)  =3$ as in \cite{dr0}.
   Then the conclusion of Theorem \ref{t-2}
   does not hold for $K=S^3$.
   In the absence of finite-dimensionality of $X$ the following
example from \cite{l0}
may serve as a source of many counter-examples: there is a compactum $X$ satisfying the following
conditions:
\begin{itemize}
\item[a.]  $\edim X > K$
   for every
   finite CW-complex $K$ with $\tilde H_*(K) \neq 0$,
\item[b.]
   $\dim_G X \leq 2$ for every abelian group $G$,
\item[c.] $\dim_G X
  \leq 1$    for every finite abelian group $G$.
\end{itemize}
   Here $\edim X > K$ means that
   $\edim X \leq K$ is false.

    With no restriction on $K$, Theorem \ref{t-2} does not hold.
  Indeed, the conclusion of Theorem \ref{t-2} is not satisfied
  if $K$ is a non-contractible  acyclic CW-complex and $X$ is
  the $2$-dimensional  disk.   Cencelj and  Dranishnikov \cite{cdr2}
  generalized  Theorem
  \ref{t-2}  for  nilpotent CW-complexes $K$
  and $X$ being a compactum (see \cite{cdr1}
for the case of $K$ with fundamental group being finitely generated).
Their work was generalized in \cite{CDMV} to $X$ being metric.
\par
The real projective plane
  $\rp^2$ is the simplest CW-complex not covered by
   Cencelj-Dranishnikov's result. Thus we arrive at the following
   well-known open problem in Extension Theory.

   \begin{Problem}
   \label{prob}
   Let $X$ be a finite dimensional compactum.
   Does $\dim_{\Z_2} (X) \leq 1$ imply $\edim X \leq \rp^2$?
   \end{Problem}
   A partial answer to \ref{prob} was given by the authors in \cite{DydLev}:

 \begin{Thm}
 \label  {t1} Let $X$ be a compactum of dimension at most three.
If $\dim_{\Z_2}(X) \leq 1$, then $\edim X \leq \rp^2$.
 \end{Thm}

 This paper is devoted to solving \ref{prob} completely. 
In view of that it is of interest to address the following question:
   
  \begin{Problem}
  \label{q1} Let $X$ be a compactum of finite dimension.
 Does $\dim_{\Z_p}(X) \leq 1$ imply  $\edim(X) \leq M(\Z_p,1)$
 for any Moore complex $M(\Z_p,1)$ with $\pi_1(M(\Z_p,1))=\Z_p$?
 \end{Problem}
 
 The big picture of our solution to \ref{prob} is as follows:
 \par Given $f\colon A\to\rp^2$ we know it extends to
 $\tilde f\colon X\to \rp^{n+2}$ for some $n\ge 1$. Our strategy is to push
 $\tilde f$ into $\rp^{n+1}$. To accomplish that we pick the perpendicular
 $\rp^n_\perp$ to $\rp^2$ in $\rp^{n+2}$ and we try to push $\tilde f$ off
 a closed tubular neighborhood $N$ of $\rp^n_\perp$.
 $\partial N$ has a natural map $\partial N\to\rp^{n+1}$
 that induces an element of the fundamental group $\pi_1(\Map(\rp^n,\rp^{n+1}),i)$
 based at inclusion $i\colon \rp^n\to\rp^{n+1}$. This is most evident
 if the circle bundle $\partial N\to\rp^n$ is trivial (for example if $n=1$).
 In that case pick $x_0\in Int(N)\setminus \rp^{n+1}$
 and observe the inclusion $\rp^{n+1}\to \rp^{n+2}\setminus\{x_0\}$ is a homotopy equivalence.
 Therefore the inclusion $\partial N\to \rp^{n+2}\setminus\{x_0\}$
 induces a map $\rp^n\times S^1\to\rp^{n+1}$
 that gives rise to an element of $\pi_1(\Map(\rp^n,\rp^{n+1}),i)$.
 That element is always notrivial
 but its square equals one. That allows to factor
 $\partial N\to\rp^{n+1}$ through a complex that can be shown to be
 an absolute extensor of $X$. Therefore $\tilde f$ can be adjusted to miss $Int(N)$.
 By induction we push $\tilde f$ into $\rp^3$ and then we use techniques
 of \cite{DydLev} to complete the proof.
  
\section{Preliminary results}
In this section we prove a few results needed for Section 
\ref{MainResults}.

Given a closed subset $L$ of $K$
by $\frac{K\times F}{L\times F}$
we mean the quotient space of $K\times F$
under the decomposition consisting of singletons
$(x,y)\notin L\times F$ and sets $\{x\}\times F$ for $x\in L$.
Notice one has a natural projection
$\pi_K\colon \frac{K\times F}{L\times F}\to K$.

\par
We use the convention $\rp^{-1}=\rp^0$ 
and consider $\rp^0$ the base-point of any $\rp^n$.
$p_n\colon S^n\to \rp^n$ is the quotient map.
We will use the same symbol to denote the quotient
map $p_n\colon B^n\to\rp^n$ from the unit $n$-ball $B^n$ onto $\rp^n$.
\par

Given any $u\colon \rp^n\times [0,1]\to \rp^{n+1}$
such that both $u_0$ and $u_1$ are inclusions
and $u\vert \rp^0\times I$ determines a homotopically trivial loop
in $\rp^{n+1}$, one can homotop $u$ relatively to $\rp^n\times\partial I$
to get $v\colon \rp^n\times [0,1]\to \rp^{n+1}$
such that $v\vert \rp^{n-1}\times I$ is the projection onto
$\rp^{n-1}$ followed by the inclusion
$\rp^{n-1}\to\rp^{n+1}$ (see \ref{NormalFormLemma}).
Let us call such $v$ a {\it normal form} of $u$
and $v$ is said to be in {\it normal form}.

Similarly, given $u\colon \rp^n\times\rp^1\to \rp^{n+1}$
such that $u\vert \rp^n\times\rp^0$ is homotopic to inclusion
$\rp^n\to\rp^{n+1}$ and $u\vert \rp^0\times \rp^1$
being homotopically trivial, one can homotop $u$
relatively to $\rp^0\times\rp^0$ to $v\colon \rp^n\times\rp^1\to \rp^{n+1}$
such that $v\vert\rp^n\times\rp^0$ is the inclusion
and $v\vert\rp^{n-1}\times\rp^1$ is the projection onto
the first coordinate. Again, such $v$ will be called a {\it normal form}
of $u$ and $v$ will be called to be in {\it normal form}.

\begin{Prop}\label{TubNConnection}
Consider a closed tubular neighborhood $N$
of $\rp^n$ in $\rp^{n+2}$ for some $n\ge 1$ and pick $x_0\in Int(N)\setminus \rp^{n+1}$.
\begin{itemize}
\item[a.] The inclusion $\rp^{n+1}\to \rp^{n+2}\setminus\{x_0\}$
is a homotopy equivalence.
\item[b.] The inclusion $\partial N\to \rp^{n+2}\setminus\{x_0\}$
factors up to homotopy
as $\partial N\to \frac{\rp^n\times\rp^1}{\rp^{n-2}\times\rp^1}\to \rp^{n+2}\setminus\{x_0\}$.
\end{itemize}
\end{Prop}
\dokaz
a). This is obvious from the representation of $\rp^{n+2}$
as the quotient of the $(n+2)$-ball $B^{n+2}$: removing a point in the interior
of $B^{n+2}$ allows for a deformation retraction to its boundary.
\par b).
Let $\pi\colon \partial N\to \rp^n$ be the projection
and let $r\colon \partial N\to \rp^{n+1}$ be the restriction
of a deformation retraction $\rp^{n+2}\setminus\{x_0\}\to\rp^{n+1}$.
Notice $r$ restricted to fibers of $\pi$ is null-homotopic
(fibers/circles are null-homotopic in $\rp^{n+2}$ and the inclusion
$\rp^{n+1}\to\rp^{n+2}$ is an isomorphism on the fundamental groups)
On $M=\pi^{-1}(\rp^{n-2})$ the map $r\vert_M$ represents an element
of $H^1(M;\Z_2)$ as $\dim(M)=n-1$, so by \ref{DetectingMapsToKGOne} it can be factored
up to homotopy through $\rp^{n-2}$. Since $\pi$ is trivial over $\rp^n\setminus\rp^{n-2}$
(see \ref{CircleBundleIsTrivial}), $r$ factors up to homotopy
as $s\circ p$, where $p\colon\partial N\to \frac{\rp^n\times\rp^1}{\rp^{n-2}\times\rp^1}$
sends fibers to fibers and $s\colon \frac{\rp^n\times\rp^1}{\rp^{n-2}\times\rp^1}\to\rp^{n+1}$.
\enddokaz

Suppose $v\colon \rp^n\times [0,1]\to \rp^{n+1}$ is in normal form for some $n\ge 1$.
Let $B^n$ be the upper hemisphere of $S^n\subset S^{n+1}$.
Consider $p_n\times id\colon B^n\times I\to \rp^n\times I$
and let $\tilde v\colon B^n\times I\to S^{n+1}$
be the lift
$$
\CD
B^n\times I @>{\tilde v}>> S^{n+1}  \\
@VV{p_n\times id}V  @VV p_{n+1} V\\
\rp^n\times I @>{v}>> \rp^{n+1}\\
\endCD
$$
of $v\circ p_n\times id$ so that $\tilde v_0$ is the inclusion $B^n\to S^{n+1}$.
Notice $\tilde v_1$ is also the inclusion and $\tilde v\vert S^{n-1}\times I$
is the projection. Thus $\tilde v(\partial (B^n\times I))\subset B^n$
and we can talk about the degree $\deg(\tilde v)$
of the induced map $(B^n\times I)/(\partial (B^n\times I))\to S^{n+1}/B^n$.

\begin{Lem}\label{BasicLemmaEvenDeg}
If $\deg(\tilde v)$ is even, then $v$ is homotopic
rel. $\rp^n\times\partial I\cup\rp^{n-2}\times I$ to the
projection $\rp^n\times I\to\rp^n$ followed by inclusion.
\end{Lem}
\dokaz
Suppose $\deg(\tilde v)=2k$.
Express $S^{n-1}$ as the union $H_+\cup H_-$ of upper hemisphere
and lower hemisphere.
Our plan is to find $G\colon I\times B^n\times I\to S^{n+1}$
so that $p_{n+1}\circ G$ will define a desired homotopy on $\rp^n\times I$.
That is accomplished by defining $G$ on $\partial(I\times B^n\times I)$
so that its degree is trivial which allows for extension of $G$ over
the entire $I\times B^n\times I$.
Since $G_0=\tilde v$ and $G_1$ are known as well as $G\vert I\times S^{n-2}\times I$,
the only two missing parts are $I\times H_+\times I$ and $I\times H_-\times I$.
Define $G$ on $I\times H_+\times I$ so that it extends $G\vert I\times S^{n-2}\times I$
and its degree is $-k$. For $x\in H_-$ we put $G(s,x,t)=-G(s,-x,t)$.
That way the degree of $G\vert I\times H_-\times I$ is also $-k$
(there are even number of negative signs in total when one does
antipodal maps on $S^{n-1}$ and $S^{n+1}$).
Thus the total degree of $G\vert \partial(I\times B^n\times I)$ is $0$
and it can be extended over $I\times B^n\times I$.
Since $x\in S^{n-1}$ implies $G(s,-x,t)=-G(s,x,t)$, $G$
indeed induces the desired homotopy between $v$
and the projection $\rp^n\times I\to \rp^{n+1}$.
\enddokaz

\begin{Prop}\label{FactThroughRP2}
Suppose $n\ge 1$.
If $\pi_1(s)$ of
$s\colon \frac{\rp^n\times\rp^1}{\rp^{n-2}\times\rp^1}\to\rp^{n+1}$
 is not trivial, then $s$ extends over
$\frac{\rp^n\times\rp^2}{\rp^{n-2}\times\rp^2}$.
\end{Prop}
\dokaz
First consider the case of the composition $u$ of
$$\rp^n\times\rp^1\overset{proj}\to \frac{\rp^n\times\rp^1}{\rp^{n-2}\times\rp^1}
\overset{s}\to\rp^{n+1}$$ being in normal form.
The composition $v$ of
$\rp^n\times S^1\overset{id\times p_1}\to \rp^n\times\rp^1\overset{u}\to\rp^{n+1}$
is of even degree, so it extends over $\rp^n\times B^2$,
$B^2$ being the upper hemisphere of $S^2$.
Let $G\colon \rp^n\times B^2\to\rp^{n+1}$ be such extension.
\ref{BasicLemmaEvenDeg} guarantees that on each slice $\{x\}\times B^2$,
$x\in \rp^{n-2}$, the map $G$ is constant.
Therefore $G$ induces $H\colon \frac{\rp^n\times\rp^2}{\rp^{n-2}\times\rp^2}\to\rp^{n+1}$
which extends $s$.
\par Consider the general case of arbitrary $s$. 
Let us show the restriction $w$ of $s$ to
$K=\frac{\rp^{n-1}\times\rp^1}{\rp^{n-2}\times\rp^1}$ is homotopic
to the projection. It is clearly so for $n\leq 1$, so assume $n\ge 2$.
Since the dimension of $K$ is $n$, $w$ represents
an element of $H^1(K;\Z_2)=[K,\rp^\infty]=[K,\rp^{n+1}]$.
Observe $s\vert \rp^{n-1}\times\rp^0$ is homotopic to the inclusion
as it is not trivial on the fundamental group (see \ref{RelRPLemma}).
Assume $s\vert \rp^{n-1}\times\rp^0$ is the inclusion.
Put $M=\rp^{n-1}\times\rp^0\subset K$
and notice $K/M$ is simply connected, so $H^1(K,M;\Z_2)=0$.
From the exact sequence $0=H^1(K,M;\Z_2)\to H^1(K;\Z_2)\to H^1(M;\Z_2)=\Z_2$
one gets $s$ is homotopic to the projection
onto the first coordinate as both restrict to the same map on $M$.

\par
So assume $s$ equals
the projection to $\rp^{n+1}$ on $\frac{\rp^{n-1}\times\rp^1}{\rp^{n-2}\times\rp^1}$.
Now $s\vert \rp^n\times\rp^0$ is homotopic to the inclusion
rel.$\rp^{n-1}\times\rp^0$ (see \ref{RelRPLemma}), so we may assume $s$
restricted to $\rp^n\times\rp^0\cup \frac{\rp^{n-1}\times\rp^1}{\rp^{n-2}\times\rp^1}$ is the projection
onto the first coordinate. That is equivalent to
the composition $u$ of
$$\rp^n\times\rp^1\overset{proj}\to \frac{\rp^n\times\rp^1}{\rp^{n-2}\times\rp^1}
\overset{s}\to\rp^{n+1}$$ being in normal form.
\enddokaz

The next result is needed to address the double surgery we do on $\rp^3$.
The first surgery yields the second modification $M_2$ of \cite{DydLev}
(see \ref{SecondModDef}). The second surgery is done on $M_2$
in a similar manner (see the proof of \ref{SecondExtLem}).
\begin{Lem}\label{FirstExtLem}
Suppose $X$ is a metric space such that $X\tau\Sigma(\rp^2)$.
 Let $r_i\colon S^1\times \rp^2\to\rp^2$, $i=1,2$,
 be maps such that $r_i\vert a\times \rp^1\sim const$
 and $r_i\vert S^1\times b$ is not null-homotopic.
  If $A$ is closed in $X$ and $f\colon A\to S^1\times \rp^1\times\rp^2$
  has the property that the composition 
  $$A\overset{f}\to S^1\times \rp^1\times\rp^2\overset{proj}\to S^1\times \rp^1$$
  extends over $X$, then the composition 
  $$A\overset{f}\to S^1\times \rp^1\times\rp^2
  \overset{id\times r_1}\to \rp^1\times\rp^2
  \overset{r_2}\to
  \rp^2$$
   extends over $X$. 
\end{Lem}
\dokaz
Notice that any composition of two maps from $\rp^2$ to $\rp^2$
that are trivial on the fundamental group, is homotopically trivial.
The reason is that each of them factors through $S^2$
and any composition $S^2\to\rp^2\S^2$ is null-homotopic
as it is of degree $0$.
Thus the composition $S^1\times \rp^1\times\rp^2
  \overset{id\times r_1}\to \rp^1\times\rp^2
  \overset{r_2}\to
  \rp^2$ is homotopically trivial on fibers of projection
  $S^1\times \rp^1\times\rp^2\overset{proj}\to S^1\times \rp^1$
  and we can apply Lemma \ref{ExtensionLem}.
\enddokaz

The second modification $M_2$ of $\rp^3$ was defined in \cite{DydLev} as follows:
\begin{Def}[Second modification of $\rp^3$]\label{SecondModDef}
Represent $\rp^3$ as the quotient of $B^3$ with $\rp^2$
being the image of $S^2=\partial B^3$.
Let $\rp_\perp^1$ be the image of the segment in $B^3$
perpendicular to the plane of $S^1$.
It has a closed tubular neighborhood 
that can be represented as a solid torus $T=S^1\times D$
with $D$ being a disk in $\rp^2$.
There is a radial retraction $r\colon \partial T\to\rp^2$
such that $r\vert S^1\times a$ is not null-homotopic
and $r\vert b\times\partial D$ is the inclusion and therefore
null-homotopic. It was shown in \cite{DydLev} (one can also use
\ref{FactThroughRP2}) that $r$ extends to
$r\colon S^1\times \rp^2\to \rp^2$,
where $\partial D$ is identified with $\rp^1$.
$M_2$ is obtained from $\rp^3\setminus Int(T)$ by attaching
$S^1\times\rp^2$ via the inclusion $S^1\times\partial D\to S^1\times\rp^2$.
The map $r$ can be now extended to $r\colon M_2\to\rp_2$
so that $r\vert\rp^1$ is the inclusion.
\end{Def}

\begin{Lem}\label{SecondExtLem}
Suppose $X$ is a metric space such that $X\tau\Sigma(\rp^2)$ and
 $X\tau(\rp^1\to\rp^3)$. Let $r\colon S^1\times S^1\to\rp^2$
 be a map such that $r\vert a\times S^1\sim const$
 and $r\vert S^1\times b$ is not null-homotopic.
  If $A$ is closed in $X$ and $f\colon A\to S^1\times S^1$
  has the property that $A\overset{f}\to S^1\times S^1\overset{proj_1}\to S^1$
  extends over $X$ , then $A\overset{f}\to S^1\times S^1  \overset{r}\to
  \rp^2$ extends over $X$. 
\end{Lem}
\dokaz
Extend the composition $A\overset{f}\to S^1\times S^1\overset{proj_2}\to S^1\overset{inclusion}\to\rp^3$
over $X$. Let $g\colon X\to\rp^3$ be such extension. That means we have an extension $F\colon X\to S^1\times\rp^3$
of $A\overset{f}\to S^1\times S^1\overset{inclusion}\rightarrow S^1\times\rp^3$.
Let $T=S^1\times D$ be the solid torus in $\rp^3$ arising from the
second modification $M_2$ of $\rp^3$ (see \ref{SecondModDef})
and let $r_M\colon S^1\times\rp^2\to \rp^2$ be the map used to define $M_2$.
That map extends to $r_M\colon M_2\to\rp^2$. Our first goal is
to find an extension $X\to S^1\times M_2$ of 
$A\overset{f}\to S^1\times S^1\overset{inclusion}\rightarrow S^1\times M_2$.
Put $Y=g^{-1}(T)$ and $C=g^{-1}(\partial T)$. 
By Lemma \ref{FirstExtLem} the composition
$C\overset{f_1\times g}\to S^1\times S^1\times \partial D
\overset{id\times r_M}\to S^1\times \rp^2
\overset{r}\to\rp^2$ extends over $Y$.
On $Z=g^{-1}(\rp^3\setminus Int(T))$ we have the map
$Z\to S^1\times M_2
\overset{id\times r_M}\to S^1\times \rp_2
\overset{r}\to\rp^2$
that can be pasted with the one on $Y$ to yield the desired $X\to\rp^2$.
\enddokaz

\begin{Lem}\label{ThirdExtLem}
Suppose $X$ is a metric space such that $X\tau\Sigma(\rp^2)$.
Let $N$ be a closed tubular neighborhood of $\rp^n$ in $\rp^{n+2}$ for some
$n\ge 2$ and let $\pi\colon\partial N\to \rp^n$ be the corresponding
circle bundle. If $f\colon X\to N$ and $r\colon\partial N\to \rp^{n+1}$
is a map such that $r\vert_F$ is null-homotopic on the fiber $F$ of $\pi$
and $\pi_1(r)$ is not trivial,
then $r\circ f\vert_A\colon A=f^{-1}(\partial N)\to \rp^{n+1}$
extends over $X$.
\end{Lem}
\dokaz
On $M=\pi^{-1}(\rp^{n-2})$ the map $r\vert_M$ represents an element
of $H^1(M;\Z_2)$ as $\dim(M)=n-1$, so by \ref{DetectingMapsToKGOne} it can be factored
up to homotopy through $\rp^{n-2}$. Since $\pi$ is trivial over $\rp^n\setminus\rp^{n-2}$
(see \ref{CircleBundleIsTrivial}), $r$ factors up to homotopy
as $s\circ p$, where $p\colon\partial N\to \frac{\rp^n\times\rp^1}{\rp^{n-2}\times\rp^1}$
sends fibers to fibers and $s\colon \frac{\rp^n\times\rp^1}{\rp^{n-2}\times\rp^1}\to\rp^{n+1}$.
As $\pi_1(s)$ is not trivial, it extends over
$\frac{\rp^n\times\rp^2}{\rp^{n-2}\times\rp^2}$ by \ref{FactThroughRP2} which completes the proof
by \ref{ExtensionLem}.
\enddokaz

\begin{Lem}\label{ReductionToLowerRP2ToRPn}
Suppose $X$ is a metric space such that $X\tau\Sigma(\rp^2)$
and $n\ge 2$.
If $X\tau(\rp^2\to\rp^{n+2})$, then $X\tau(\rp^2\to\rp^{n+1})$. 
\end{Lem}
\dokaz
Pick perpendicular $\rp^n_\perp$ to $\rp^2$ in $\rp^{n+2}$.
Let $N$ be a closed tubular neighborhood
of $\rp^n_\perp$ in $\rp^{n+2}$ and let $r\colon \partial N\to M=\rp^{n+2}\setminus {a_0}$
be the the inclusion where $a\in N\setminus (\partial N\cup \rp^{n+1})$.
Notice the inclusion $\rp^{n+1}\to M$ is a homotopy equivalence.
Also notice $\pi_1(r)$ is not trivial as $\partial N$ has a loop homotopic
in $N$ to the generator of $\pi^{n+2}$, so the inclusion
$\partial N\to \rp^{n+2}$ is non-trivial on the level of fundamental groups.

\par Given $f\colon A\to\rp^2$ extend it to $F\colon X\to \rp^{n+2}$
and put $Y=F^{-1}(N)$, $C=F^{-1}(\partial N)$.
By \ref{ThirdExtLem} there is $F\colon X\to M$
agreeing with $f$ on $F^{-1}(\rp^{n+2}\setminus (N\setminus\partial N))$.
Let $D$ be the disk $\rp^2\cap N$.
Notice $X\overset{F}\to M\to M/D$ extends $A\overset{f}\to\rp^2\to\rp^2/D$.
Since the the projection $(M,\rp^2)\to (M/D,\rp^2/D)$ and 
the inclusion $(\rp^{n+1},\rp^2)\to (M,\rp^2)$
are homotopy equivalences, the proof is completed.
\enddokaz

\section{Main results}\label{MainResults}

\begin{Thm}\label{GenThmForRP2}
Suppose $X$ is a metric space.
If $X\tau(\rp^1\to\rp^2)$ and $X\tau(\rp^2\to\rp^3)$,
then $X\tau \rp^2$. 
\end{Thm}
\dokaz
Suppose $f\colon A\to \rp^2$, $A$ being closed in $X$.
Extend $f$ to $F\colon X\to\rp^3$. 
Represent $\rp^3$ as the quotient of $B^3$ with $\rp^2$
being the image of $S^2=\partial B^3$.
Remove a solid torus $T=S^1\times D$ in the interior of $B^3$.
Attach $S^1\times \rp^2$ to $B^3\setminus Int(T)$ via
identity on $S^1$ times the inclusion $\partial D\to\rp^2$,
where $\partial D$ is identified with $\rp^1$.
The resulting space $M$ retracts onto $\rp^2$ as shown in \cite{DydLev}
(it is called the {\it first modification of $\rp^3$} there).
Put $C=F^{-1}(T)$ and $D=F^{-1}(\partial T)$.
Since $C\tau(\rp^1\to\rp^2)$, $F\vert_D\colon D\to \partial T$
extends to $G\colon C\to S^1\times\rp^2$.
Replace $F$ by $G$ on $C$ and follow the retraction
of $M$ onto $\rp^2$ to get an extension
of $f$ from $X$ to $\rp^2$.
\enddokaz

\begin{Thm}\label{ThmForRP1}
Suppose $X$ is a metric space such that $X\tau\Sigma(\rp^2)$.
If $X\tau(\rp^1\to\rp^3)$, then $X\tau(\rp^1\to\rp^2)$. 
\end{Thm}
\dokaz
It suffices to show $X\tau(\rp^1\to M)$, $M$ being the second modification,
as the map $r_M\colon M\to \rp^2$ gives an extension $r_m\circ F\colon X\to \rp^2$
of any $f\colon A\to \rp^1$ where $F\colon X\to M$ extends $f$.
\par Given $f\colon A\to \rp^1$ and its extension $F\colon X\to\rp^3$
put $Y=F^{-1}(T)$, $T=S^1\times D$ being the solid torus from the definition of
$M$. Since $F^{-1}(\partial T)\to S^1\times \partial D\overset{r_M}\to \rp_2$
extends over $Y$ to $g\colon Y\to\rp^2$ by Lemma \ref{SecondExtLem}, one can paste  
$Y\to S^1\times \rp^2\subset M$ with $F$ restricted to $X\setminus F^{-1}(Int(T))$
to obtain an extension $X\to M$ of $f$.
\enddokaz

\begin{Cor}\label{MainCorollary}
Suppose $X$ is a metric space. If $X\tau\Sigma(\rp^2)$,
then the following conditions are equivalent:
\begin{itemize}
\item[1.] $X\tau \rp^2$,
\item[2.] If $X\tau(\rp^2\to\rp^3)$. 
\item[3.] If $X\tau(\rp^2\to\rp^n)$ for some $n\ge 3$. 
\end{itemize}
\end{Cor}
\dokaz
It suffices to prove 2)$\implies$1) as 3)$\implies$2) is shown in \ref{ReductionToLowerRP2ToRPn}.
By Theorem \ref{ThmForRP1} one has $X\tau(\rp^1\to\rp^2)$
which combined with Theorem \ref{GenThmForRP2} yields $X\tau \rp^2$.
\enddokaz

\begin{Thm}\label{MainTheorem}
Suppose $X$ is a metric space. If $X$ is of finite dimension,
then the following conditions are equivalent:
\begin{itemize}
\item[1.] $X\tau \rp^2$,
\item[2.] $X\tau\rp^\infty$. 
\end{itemize}
\end{Thm}
\dokaz
1)$\implies$2) follows from \ref{t-1}.
\par 2)$\implies$1). Notice
 $\Sigma(\rp^2)$ is a simply-connected Moore space
$M(\Z_2,2)$, so applying \ref{t-2} for metric spaces \cite{Dy1} gives $X\tau\Sigma(\rp^2)$.
Use \ref{MainCorollary}.
\enddokaz

\section{The fundamental group of function spaces}

Consider the inclusion $i\colon\rp^n\to\rp^{n+1}$. If $n\ge 1$,
then one has the homomorphism
 $$\pi_1(\Map(\rp^n,\rp^{n+1}),i)\to \Z_2$$
 induced by the evaluation map (at $\rp^0$) $\Map(\rp^n,\rp^{n+1})\to \rp^{n+1}$.
 The goal of this section is show how the kernel of this homomorophism is 
 related to \ref{prob}.
We show the homomorphism is an epimorphism and its kernel is $\Z_2$.

\begin{Prop}\label{TubNConnection2}
Consider a closed tubular neighborhood $N$
of $\rp^n$ in $\rp^{n+2}$ for some $n\ge 1$ and pick $x_0\in Int(N)\setminus \rp^{n+1}$.
\begin{itemize}
\item[a.]  Any map $u\colon \frac{\rp^n\times\rp^1}{\rp^{n-2}\times\rp^1}\to \rp^{n+1}$
induces an element of 
\par\noindent
 $\pi_1(\Map(\rp^n,\rp^{n+1}),i)$
that is trivial if and only if $u$ is homotopic to the projection onto
the first coordinate.
\item[b.] The element of $\pi_1(\Map(\rp^n,\rp^{n+1}),i)$
induced by $\partial N$ is non-trivial and belongs to the kernel
of  $\pi_1(\Map(\rp^n,\rp^{n+1}),i)\to \Z_2$.
\end{itemize}
\end{Prop}
\dokaz
a).  Consider the projection
$p\colon \rp^n\times\rp^1\to \frac{\rp^n\times\rp^1}{\rp^{n-2}\times\rp^1}$.
Any $u\colon \frac{\rp^n\times\rp^1}{\rp^{n-2}\times\rp^1}\to \rp^{n+1}$
gives rise to the composition $u\circ p\colon \rp^n\times\rp^1\to \rp^{n+1}$
and that
induces an element of 
 $\pi_1(\Map(\rp^n,\rp^{n+1}),i)$.
Using \ref{BigHomotopyLem} one can see that element is trivial if and only if $u$ is homotopic to the projection onto
the first coordinate.
\par
b). Clearly the element of $\pi_1(\Map(\rp^n,\rp^{n+1}),i)$
induced by $\partial N$ belongs to the kernel
of  $\pi_1(\Map(\rp^n,\rp^{n+1}),i)\to \Z_2$.
If it was trivial, then the inclusion $\partial N\to \rp^{n+1}\setminus\{x_0\}$
could be factored up to homotopy (use \ref{BigHomotopyLem})
through the projection $\partial N\to \rp^n$ and would be extendable over $N$.
That would lead to a map $\rp^{n+2}\to\rp^{n+1}$ that is non-trivial
on the fundamental group contradicting Borsuk-Ulam Theorem
(after lifting to the covering spheres).
\enddokaz

\begin{Thm}\label{KernelOfHomTheorem}
Consider the inclusion $i\colon\rp^n\to\rp^{n+1}$. If $n\ge 1$,
then the homomorphism
 $$\pi_1(\Map(\rp^n,\rp^{n+1}),i)\to \Z_2$$
 induced by the evaluation map (at $\rp^0$) $\Map(\rp^n,\rp^{n+1})\to \rp^{n+1}$ is 
 an epimorphism and its kernel is $\Z_2$.
\end{Thm}
\dokaz
As $\rp^\infty$ is an Eilenberg-MacLane space $K(\Z_2,1)$
there is a pointed cellular map $m\colon \rp^\infty\times\rp^\infty\to\rp^\infty$
corresponding to addition $\Z_2\times\Z_2\to\Z_2$ on the level of fundamental
groups. Since $m(\rp^n\times\rp^1)\subset\rp^{n+1}$,
the restriction $m\vert (\rp^n\times\rp^1)$ of $m$ induces a loop
in $(\Map(\rp^n,\rp^{n+1}),i)$ whose evaluation is a non-trivial loop in
$\rp^{n+1}$. That proves the homomorphism
 $\pi_1(\Map(\rp^n,\rp^{n+1}),i)\to \Z_2$
 induced by the evaluation map (at $\rp^0$) $\Map(\rp^n,\rp^{n+1})\to \rp^{n+1}$ is 
 an epimorphism.
 \par
Every loop $S^1\to \Map(\rp^n,\rp^{n+1})$ based at inclusion
can be converted to a map $\alpha\colon \rp^n\times S^1\to \rp^{n+1}$
such that $\alpha\vert (\rp^n\times 1)$ is the inclusion.
$\alpha$ belongs to the  kernel of the homomorphism
 $$\pi_1(\Map(\rp^n,\rp^{n+1}),i)\to \Z_2$$
 induced by the evaluation map if and only if
 $\alpha\vert (\rp^0\times S^1)$ is null-homotopic.
 By  \ref{NormalFormLemma} we may assume that $\alpha$ is in normal form
 and \ref{BasicLemmaEvenDeg} says that there are at most two homotopy classes
of such maps depending on whether the degree of $\alpha$
is even or not (notice the degree is additive under loop concatenation for loops in normal form).
A geometrical way to detect non-trivial element of the kernel follows from
Part d) of \ref{TubNConnection}.
Here is another one:
 pick $v$ in normal form so that the degree of its lift $\tilde v$
 $$
\CD
B^n\times I @>{\tilde v}>> S^{n+1}  \\
@VV{p_n\times id}V  @VV p_{n+1} V\\
\rp^n\times I @>{v}>> \rp^{n+1}\\
\endCD
$$
is $1$.
 If $\alpha$ induces the trivial element of $\pi_1(\Map(\rp^n,\rp^{n+1}),i)$,
then there is a cellular homotopy (rel.$\partial I$) $H\colon I\times \rp^n\times I\to \rp^{n+1}$
joining $v$ and the projection $p\colon \rp^n\times I\to\rp^{n+1}$
onto the first coordinate.
Consider the lift
$\tilde H\colon I\times B^n\times I\to S^{n+1}$ of $H$.
$\tilde H_0\colon \{0\}\times B^n\times I\to S^{n+1}$ equals $\tilde v$,
$\tilde H_1\colon \{1\}\times B^n\times I\to S^{n+1}$ equals
the projection onto $B^n$, and both $\tilde H\vert (I\times B^n\times \{0\}$
and $\tilde H\vert (I\times B^n\times \{1\}$ are projections onto $B^n$.
Therefore we need only to look at $\tilde H\vert (I\times S^{n-1}\times I$
to compare its the degree to that of $\tilde v$.
Let $H_+$ be the upper hemisphere of $S^{n-1}$ and let
$H_-$ be  the lower hemisphere of $S^{n-1}$.
If $x\in H_+$, then $\tilde H(s,-x,t)=\pm\tilde H(s,x,t)$ for all $s,t\in I$.
Therefore the cumulative degree coming from $\tilde H\vert (I\times H_+\times I)$
and $\tilde H\vert (I\times H_-\times I)$ is even, a contradiction.
 \enddokaz
 
 In view of \ref{KernelOfHomTheorem} there are only two possibilities
 for the group (given $n\ge 1$) $\pi_1(\Map(\rp^n,\rp^{n+1}),i)$: $\Z_4$ or $\Z_2\oplus\Z_2$.
 We do not know which case holds for a particular $n$,
 only for some initial values of $n$.
   \begin{Prop}\label{CharOfZ2Z2VsZ4}
If $n\ge 1$, then the following conditions are equivalent:
\begin{itemize}
\item[a.] $\pi_1(\Map(\rp^n,\rp^{n+1}),i)=\Z_2\oplus\Z_2$.
\item[b.] There
is $u\colon \rp^n\times\rp^2\to \rp^{n+1}$
such that both $u\vert (\rp^n\times\rp^0)$
and $u\vert (\rp^0\times\rp^2)$ are inclusions.
\end{itemize}
 \end{Prop}
 \dokaz
 a)$\implies$b). Pick 
 $u\colon \rp^n\times \rp^1\to\rp^{n+1}$ not belonging to the kernel
 of evaluation-induced epimorphism $\pi_1(\Map(\rp^n,\rp^{n+1}),i)\to\Z_2$.
 We may assume $u\vert\rp^0\times\rp^1$ is the inclusion.
 Since $u$ is of degree $2$, it extends over
$\rp^n\times\rp^2$.
Since any two maps $\rp^2\to\rp^{n+1}$ being inclusions on $\rp^1$
are homotopic (see \ref{RelRPLemma}), we may assume the extension
equals inclusion on $\rp^0\times\rp^2$.
\par
b)$\implies$a). Any such $u$ induces an element of order $2$
not belonging to the kernel
 of evaluation-induced epimorphism $\pi_1(\Map(\rp^n,\rp^{n+1}),i)\to\Z_2$.
 That means $\pi_1(\Map(\rp^n,\rp^{n+1}),i)$ cannot be $\Z_4$.
 \enddokaz
 
 \begin{Cor}
$\pi_1(\Map(\rp^n,\rp^{n+1}),i)=\Z_2\oplus\Z_2$
for $n=2,3,6,7$.
 \end{Cor}
 \dokaz
The quaternionic multiplication $S^3\times S^3$
induces $u\colon \rp^3\times\rp^3\to\rp^3$
such that both $u\vert (\rp^3\times\rp^0)$
and $u\vert (\rp^0\times\rp^3)$ are inclusions.
This takes care of $n=2,3$.
\par Using octonions (the Cayley numbers) and the multiplication
$S^7\times S^7\to S^7$ one induces
$u\colon\rp^7\times\rp^7\to \rp^7$
which handles $n=6,7$.
 \enddokaz

 \begin{Cor}
$\pi_1(\Map(\rp^1,\rp^{2}),i)=\Z_4$.
 \end{Cor}
 \dokaz
 If $\pi_1(\Map(\rp^1,\rp^{2}),i)=\Z_2\oplus\Z_2$, then
 \ref{CharOfZ2Z2VsZ4} implies
 existence of $u\colon \rp^1\times\rp^2\to \rp^2$
 such that $u\vert(\rp^1\times\rp^0)$ and $u\vert(\rp^0\times\rp^2)$
 are inclusions. Converting $u$ to $u\colon I\times\rp^2\to\rp^2$
 and lifting it to $\tilde u\colon I\times S^2\to S^2$
 produces a homotopy joining identity on $S^2$ with the antipodal map,
 a contradiction.
 \enddokaz

\section{Appendix A}

In this section we prove results on fibrations and bundles
that are necessary for the paper.

\begin{Prop}\label{FactoringThroughFunnyStuff}
Suppose $p\colon E\to B$ is a bundle with compact fiber $F$
that is trivial over $B\setminus C$, $C$ being closed in $B$.
If $f\colon E\to K$ is a map such that $f(p^{-1}(c))$ is a point
for each $c$, then 
there are maps $p'\colon E\to\frac{B\times F}{C\times F}$
and $g\colon \frac{B\times F}{C\times F}\to K$
such that

$$
\CD
E @>{ p'}>> \frac{B\times F}{C\times F}  \\
@VV{p}V  @VV \pi_B V\\
B @>{id_B}>> B\\
\endCD
$$
is commutative, $f=g\circ p'$,
and $p'$ is the isomorphism of bundles over $B\setminus C$.
\end{Prop}
\dokaz Choose a trivialization $h\colon (B\setminus C)\times F\to p^{-1}(B\setminus C)$
of $p$ over $B\setminus C$.
$p'$ equals $h^{-1}$ on $p^{-1}(B\setminus C)$
and sends each $p^{-1}(c)$, $c\in C$, to the corresponding point of 
$\frac{B\times F}{C\times F}$. 
The only item to check is the continuity of $p'$ at points in $p^{-1}(C)$.
That follows from the fact each neighborhood of $c\times F$ in $B\times F$
contains $U\times F$ for some neighborhood $U$ of $c$ in $B$.
\enddokaz

\begin{Thm}[\cite{DD}]\label{MainLiftingThm}
Suppose $B$ is a simplicial complex with the weak topology
and $p\colon E\to B$ is a map.
If $X$ is a metric space such that $p^{-1}(\Delta)$ is an absolute extensor of $X$
for all simplices $\Delta$ of $B$, then for any commutative diagram
$$
\CD
A @>{\tilde f}>> E  \\
@VV{i}V  @VV p V\\
X @>{g}>> B\\
\endCD
$$
where $A$ is closed in $X$, there is an extension
$F\colon X\to E$ of $f$ such that $F(x)$ and $g(x)$ belong
to the same simplex $\Delta_x$ of $B$ for all $x\in X$
(in particular, $g$ is homotopic to $p\circ F$ relatively the set of points
on which they coincide).
\end{Thm}

\begin{Cor}\label{ACorollary}
Let $n\ge 0$.
Suppose $B$ is a simplicial complex with the weak topology
and $p\colon E\to B$ is a map such that
 $p^{-1}(\Delta)$ is $n$-connected CW complex for all simplices $\Delta$ of $B$,
 then $\pi_k(p)\colon \pi_k(E,e_0)\to\pi_k(B,p(e_0))$
 is a isomorphism for $k \leq n$ and an epimorphism for $k=n+1$.
\end{Cor}
\dokaz
Notice $p^{-1}(\Delta)$ is an absolute extensor
of any $(n+1)$-dimensional space $X$ for all simplices $\Delta$ of $B$.
Applying \ref{MainLiftingThm} for $X=S^k$ one gets
$\pi_k(p)\colon \pi_k(E,e_0)\to\pi_k(B,p(e_0))$
being an epimorphism for $k\leq n+1$.
Applying \ref{MainLiftingThm} for $X=B^{k+1}$ and $A=\partial B^{k+1}$ one gets
$\pi_k(p)\colon \pi_k(E,e_0)\to\pi_k(B,p(e_0))$
being a monomorphism for $k \leq n$.
\enddokaz

\begin{Cor}\label{DetectingMapsToKGOne}
Suppose $B$ is a finite simplicial complex.
If $p\colon E\to B$ is a circle bundle with fiber $F$,
then $p$ induces bijection of pointed homotopy classes
$[E/F,K(G,1)]$ and $[B,K(G,1)]$ for any group $G$.
In particular, if $G$ is Abelian,
$p$ induces bijection of free homotopy classes
$[E/F,K(G,1)]$ and $[B,K(G,1)]$
\end{Cor}
\dokaz
Pick a maximal tree $T$ in $B^{(1)}$. Notice $p^{-1}(T)/F$ is contractible
and consider $E'$ obtained from $E$ by contracting each fiber over every point
of $T$. The natural map $E/F\to E'$ is a homotopy equivalence,
and the natural projection $q\colon E'\to B$ has the property that
$q^{-1}(\Delta)$ is connected and simply connected for all
simplices $\Delta$ of $B$. Apply \ref{ACorollary}
to conclude $E/F\to B$ induces isomorphism of fundamental groups.
Since the pointed homotopy classes $[L,K(G,1)]$ are identical with
homomorphisms from $\pi_1(L)$ to $G$ if $L$ is a CW complex,
we are done.
\enddokaz

\begin{Lem}\label{ExtensionLem}
Suppose $L\ne\emptyset $ is a subcomplex
of a connected
simplicial complex $K$ and $X$ is a metric space such that
$X\tau\Sigma(F)$ for some CW complex $F$. Let $u\colon \frac{K\times F}{L\times F}\to M$ be a map
to a CW complex.
Given a closed subset $A$ of $X$ and a map
$f\colon A\to \frac{K\times F}{L\times F}$ such that $\pi_K\circ f$ extends over $X$,
$u\circ f$ extends over $X$.
\end{Lem}
\dokaz
First, triangulate $K$ so that $L$ is a full subcomplex of $K$.
Second, since $u$ is null-homotopic on every slice $\{x\}\times F$, $x\in K$,
add to $L$ all the vertices of $K$. This way $L$ minus isolated points
is a full subcomplex of $K$.
In that case it is sufficient to consider $M=\frac{K\times F}{L\times F}$
and $u$ being the identity map.
For $K$ being a simplex this amounts to showing $M$
is an absolute extensor of $X$. Once this is shown, the general case
follows as $\pi^{-1}_K(\Delta)$ is an absolute extensor of $X$ for all
simplices $\Delta$ of $K$.
\par First, consider the case $K=v\ast\Delta$ and $L=\Delta\cup\{v\}$.
If $\dim(\Delta)=0$, then $\frac{K\times F}{L\times F}$ is simply
the suspension of $F$. Higher dimensional cases
reduce to lower ones by using a deformation retraction
of $\Delta$ to its face. That deformation retraction
extends to $\frac{K\times F}{L\times F}$.
Thus $\frac{K\times F}{L\times F}$ is homotopy equivalent to $\Sigma(F)$
if $K=v\ast\Delta$ and $L=\Delta\cup\{v\}$.
If $K=\Delta$, consider the barycenter $v$ of $\Delta$
and let $K'$ be the triangulation of $K$ in the form of $v\ast \sigma$
for all proper faces $\sigma$ of $\Delta$. If we put $L'=L\cup\{v\}$,
then, by the previous case, the projection
$\pi_{K'}\colon \frac{K'\times F}{L'\times F}\to  K'$
has the property that point-inverses of simplices are absolute extensors of $X$.
Thus $E'=\frac{K'\times F}{L'\times F}$ is an absolute extensor of $X$.
Observe that $E'=E/F$, where $E=\frac{K\times F}{L\times F}$ and $F$
is contractible in $F$. Since $E'$ is homotopy equivalent to $E\cup Cone(F)$
and one can retract $E\cup Cone(F)$ onto $E$, the proof is completed.
\enddokaz

\section{Appendix B}
In this section we prove a few technical results on projective spaces
that are needed and which are of general nature (not related to Extension Theory).
\begin{Lem}\label{CircleBundleIsTrivial}
Let $N$ be a closed tubular neighborhood of $\rp^n$ in $\rp^{n+2}$
 and let $\pi\colon\partial N\to \rp^n$ be the corresponding
circle bundle. If $n\ge 2$, then $\pi$ is trivial over $\rp^n\setminus\rp^{n-2}$.
\end{Lem}
\dokaz
Represent $\rp^{n+2}$ as equivalence classes $[x_1,x_2,x_3,x_4,\ldots,x_{n+3}]$
obtained from identifying antipodal points in the unit $(n+2)$-sphere.
Under this model $\rp^n$ is the set of points $[0,0,x_3,x_4,\ldots,x_{n+3}]$
and $\rp^n\setminus\rp^{n-2}$ is the set of points $[0,0,x_3,x_4,\ldots,x_{n+3}]$
such that $x_3^2+x_4^2\ne 0$.
We will consider $N$ to be the set of points
$[x_1,x_2,x_3,x_4,\ldots,x_{n+3}]$ such that $x_1^2+x_2^2\leq \frac{1}{9}$
and the disk bundle $\pi\colon N\to \rp^n$ is given
by $\pi([x_1,x_2,x_3,x_4,\ldots,x_{n+3}])=[0,0,t\cdot x_3,t\cdot x_4,\ldots,t\cdot x_{n+3}]$,
where $t=\sqrt{\frac{1}{1-x_1^2-x_2^2}}$.
Let us show the corresponding circle bundle $\pi\colon\partial N\to \rp^n$
is trivial over $\rp^n\setminus\rp^{n-2}$ by exhibiting
maps $\phi\colon \pi^{-1}(\rp^n\setminus\rp^{n-2})\to (\rp^n\setminus\rp^{n-2})\times S^1$
and $\phi\colon  (\rp^n\setminus\rp^{n-2})\times S^1\to \pi^{-1}(\rp^n\setminus\rp^{n-2})$
that are inverse to each other.
$\phi$ is defined by the formula
$$\phi([x_1,x_2,x_3,x_4,\ldots,x_{n+3}])=([0,0,t\cdot x_3,t\cdot x_4,\ldots,t\cdot x_{n+3}], \frac{z}{|z|})$$
where $t=\sqrt{\frac{1}{1-x_1^2-x_2^2}}=\sqrt{\frac{9}{8}}$
and $z=(x_1+i\cdot x_2)\cdot (x_3+i\cdot x_4)$.
$\psi$ is defined by the formula
$$\psi([0,0,x_3,x_4,\ldots,x_{n+3}],w)=[x_1,x_2,s\cdot x_3, s\cdot x_4,\ldots, s\cdot x_{n+3}]$$
where $s=\sqrt{\frac{8}{9}}$ and $x_1+i\cdot x_2=\frac{z}{3|z|}$, $z=\frac{w}{x_3+i\cdot x_4}$.
\enddokaz

\begin{Lem}\label{RelRPLemma}
Any two maps $u,v\colon \rp^n\to \rp^{n+1}$
such that $u\vert\rp^{n-1}=v\vert\rp^{n-1}$ are
are homotopic rel.$\rp^{n-1}$.
\end{Lem}
\dokaz
Consider $p_n\colon B^n\to \rp^n$
and define $H\colon \partial (B^n\times I)\to\rp^{n+1}$
as follows:
\begin{enumerate}
\item $H\vert (B^n\times 0)=u\circ p_n$.
\item $H\vert (B^n\times 1)=v\circ p_n$.
\item $H(x,t)=u(p_n(x))=v(p_n(x))$ for all $(x,t)\in S^{n-1}\times I$.
\end{enumerate}
$H$ extends over $B^n\times I$ and induces a homotopy
rel.$\rp^{n-1}$ of $u$ and $v$.
\enddokaz

\begin{Lem}\label{NormalFormLemma}
Given a homotopy $u\colon \rp^n\times [0,1]\to \rp^{n+2}$
such that both $u_0$ and $u_1$ are inclusions
and $u\vert \rp^0\times I$ determines a homotopically trivial loop
in $\rp^{n+2}$, one can homotop $u$ relatively to $\rp^n\times\partial I$
to the projection onto
$\rp^{n}$ followed by the inclusion
$\rp^{n}\to\rp^{n+1}$.
\end{Lem}
\dokaz
It is clearly so for $n=0$ and our proof will be by induction on $n$.
Therefore we may assume $u\vert \rp^{n-1}\times I$ is the projection
onto the first coordinate.
Let $B^n$ be the upper hemisphere of $S^n\subset S^{n+2}$.
Consider $p_n\times id\colon B^n\times I\to \rp^n\times I$
and let $\tilde u\colon B^n\times I\to S^{n+2}$
be the lift
$$
\CD
B^n\times I @>{\tilde u}>> S^{n+2}  \\
@VV{p_n\times id}V  @VV p_{n+2} V\\
\rp^n\times I @>{u}>> \rp^{n+2}\\
\endCD
$$
of $u\circ (p_n\times id)$ so that $\tilde u_0$ is the inclusion $B^n\to S^{n+2}$.
Notice $\tilde u_1$ is also the inclusion and $\tilde u\vert (\partial B^n)\times I$
is the projection onto the first coordinate.
Define $H\colon 
I\times (\partial B^n)\times I\cup I\times B^n\times \{0\}\cup I\times B^n\times \{1\}\cup \{0,1\}\times B^n\times I
\cup I\times\rp^0\times I\to S^{n+2}$
as the necessary prelude of the homotopy (rel. $(\partial B^n)\times I$) joining $\tilde u$ and the projection
$B^n\times I\to S^{n+2}$ onto the first coordinate.
$H$ extends over $I\times B^n\times I$ as $S^{n+2}$ is $(n+1)$-connected
and it induces a homotopy $G$ (rel.$\rp^n\times\partial I$) from $u$ to the projection onto the first coordinate.
\enddokaz

\begin{Lem}\label{BigHomotopyLem}
Suppose $f\colon \rp^1_1\times\rp^{n-2}\times\rp^1_2\to\rp^{n+1}$
is a map for some $n\ge 2$ such that
$f\vert (\rp^0_1\times\rp^{n-2}\times\rp^1_2)$ is the projection onto
the second coordinate.
If both $f\vert(\rp^1_1\times\rp^0\times\rp^0_2)$
and $f\vert(\rp^0_1\times\rp^{0}\times\rp^1_2)$ are null-homotopic,
then $f$ is homotopic rel.$\rp^0_1\times\rp^{n-2}\times\rp^1_2$
to the projection onto the second coordinate.
\end{Lem}
\dokaz The proof is by induction on $n$.
For $n=2$ there is a lift $\tilde f\colon \rp^1_1\times\rp^{0}\times\rp^1_2\to S^{3}$
of $f$ and we want to construct a homotopy
$H\colon \rp^1_1\times\rp^{0}\times\rp^1_2\times I\to S^{3}$
rel.$\rp^0_1\times\rp^{0}\times\rp^1_2$
starting from $\tilde f$ and ending at the constant map to the point
$\tilde f(\rp^0_1\times\rp^{0}\times\rp^0_2)$.
Such $H$ exists as $S^3$ is an absolute extensor of
$3$-dimensional $\rp^1_1\times\rp^{0}\times\rp^1_2\times I$.
\par
For $n\ge 3$ consider the lift
$\tilde f\colon \rp^1_1\times B^{n-2}\times\rp^1_2\to S^{n+1}$
of the composition $f\circ (id\times p_{n-2}\times id)$
that is the projection onto the second coordinate
on $\rp^0_1\times B^{n-2}\times\rp^1_2$.
We want to construct a homotopy
$H\colon \rp^1_1\times B^{n-2}\times\rp^1_2\times I\to S^{n+1}$
rel.$\rp^0_1\times B^{n-2}\times\rp^1_2$
starting from $\tilde f$ and ending at the projection onto the second coordinate.
We already have such homotopy on
$\rp^1_1\times S^{n-3}\times\rp^1_2\times I$ by inductional assumption,
so we will use it in order for the resulting $H$ to descend
to a homotopy on $\rp^1_1\times \rp^{n-2}\times\rp^1_2\times I$.
We can piece together those maps and then extend them
over $\rp^1_1\times B^{n-2}\times\rp^1_2\times I$.
Such extension exists as $S^{n+1}$ is an absolute extensor of
$(n+1)$-dimensional $\rp^1_1\times B^{n-2}\times\rp^1_2\times I$.
\enddokaz

\end{document}